\documentclass{article}

\usepackage[english]{babel}

\usepackage[letterpaper,top=2cm,bottom=2cm,left=3cm,right=3cm,marginparwidth=1.75cm]{geometry}

\usepackage{amsmath,amsthm}
\usepackage{graphicx}
\usepackage{paralist}
\usepackage{todonotes}
\usepackage{tikz,tikz-cd}

\usepackage{stex}
\usepackage{booktabs}
\newcommand{\raa}[1]{\renewcommand{\arraystretch}{#1}}
\usepackage{enumitem}
\makeatletter
\def\namedlabel#1#2{\begingroup
    #2%
    \def\@currentlabel{#2}%
    \phantomsection\label{#1}\endgroup
}
\makeatother
\usepackage[colorlinks=true, allcolors=blue]{hyperref}

\title{Fraudulent Publishing in the Mathematical Sciences}

\author{Ilka Agricola, Lynn Heller, Wil Schilders, Moritz Schubotz, Peter Taylor, Luis Vega}

\date{September 2025}

\begin{document}
\maketitle

\section{Introduction}
This report is the first of two publications of a joint Working Group of the International Mathematical Union (IMU) and the International Council of Industrial and Applied Mathematics (ICIAM). In it, we shall analyze the current state of publishing in the mathematical sciences and  explain the resulting problems.
Our second publication will offer concrete recommendations, guidelines, and best practices for researchers, policymakers, and evaluators of mathematical research \cite{WG-2}. It will explain how to detect and counteract attempts to game bibliometric measures, empowering the community to reclaim control over research evaluation and drive necessary change.
%

\medskip
As students, teachers, authors, reviewers, editors and committee members, we contribute to the creation, publication, consumption and evaluation of research outputs. In the current academic environment, the assessment of research quality is heavily influenced by bibliometric analysis. As a result, there is a strong incentive for researchers to optimize their bibliometric data. Unfortunately there is a growing number who are using fraudulent methods to do this. New (and to many unknown) threats such as {\it paper mills}, {\it predatory journals} and {\it citation cartels} have been invading the `ecosystem' of scientific publishing and putting the integrity of our profession at risk. A detailed glossary is provided in Appendix A, with an extensive
list of references pointing to further reading that supports the claims made in this article. As a main source, we recommend the excellent book \cite{BL}.
It is worth noting that most of the articles in the reference list have appeared in the last five years, reflecting the urgency of the problem.
It will become clear that two driving forces behind the problems we address are
the advent of Article Processing Charges (APCs) for Open Access as a business model for scientific publishing and the pervasive demand for quantitative research output assessment. Also intricately involved are the companies that collect citation data and use it to produce measures such as the journal impact factor, the h-index, and the lists of highly-cited researchers (HCRs) that are derived from this data. We will argue that products such as these are having detrimental effects on the publication enterprise.

\medskip
This is the digital arXiv version of this article with
complete clickable references. The print version with a
selection of the most important references appeared in
the October 2025 issue of the Notices of the AMS \cite{AHSSTV}.

\section{Bibliometrics in the mathematical sciences}

\begin{quote}
\emph{When a measure becomes a target, it ceases to be a good measure.}
\end{quote}
\hfill  Charles Goodhart, British economist \cite{Goodhart2}\\

In 2009, Adler, Ewing and Taylor were engaged respectively by the Institute of Mathematical Statistics (IMS), the International Mathematical Union (IMU) and the International Council of Industrial and Applied Mathematics (ICIAM) to look into the effects that quantitative assessment of research was having on their disciplines. The charge that they were asked to deliver on was: 

\begin{quote}
\it{The drive towards more transparency and accountability in the academic world has created a “culture of numbers” in which institutions and individuals believe that fair decisions can be reached by algorithmic evaluation of some statistical data; unable to measure quality (the ultimate goal), decision-makers replace quality by numbers that they can measure. This trend calls for comment from those who professionally “deal with numbers" -
mathematicians and statisticians.}
\end{quote}
The result was the report \cite{AET}, which was published in a number of outlets. 
While the report remains highly-relevant, a new phenomenon has emerged that merits a similar discussion, the widespread gaming of quantitative metrics. To provide background, we shall in this section briefly recall what \cite{AET} and other papers \cite{CareyCowlingTaylor,HT, DMV} say about the publishing culture of the mathematical sciences and pose some questions to motivate readers' thinking. 

\subsection{How mathematicians publish}
Mathematicians typically publish smaller numbers of papers, with fewer coauthors and citations than most other scientists. This renders meaningless any comparison of mathematicians with other scientists according to bibliometric measures. Our impression is that this point is largely understood by the proponents of bibliometric analysis: they are usually willing to concede that the disciplinary context matters. 
What is less understood is that there can be very different citation cultures between different sub-disciplines. Furthermore, the small absolute numbers of citations that are a result of the mathematical science referencing culture raise the question of whether it is possible to discern any signal of `quality' from the background noise.
For us, the difficult questions revolve around how relevant bibliometrics is to evaluating research within the mathematical sciences. For example: 
\begin{enumerate}

    \item[Q 1:] Does a relatively high number of citations tell us that a paper in the mathematical sciences is better than another paper?
    \item[Q 2:] Does the impact factor of a journal in the mathematical sciences reflect its quality? 
    \item[Q 3:] Does a relatively high number of papers tell us that a researcher in the mathematical sciences is better than another researcher?
    \item[Q 4:] Should citation analysis over a period of time be used to evaluate the careers of mathematical scientists?
    \item[Q 5:] Should aggregated citation data be used to rank institutions in the mathematical sciences?
    \item[Q 6:] Does the existence of bibliometric measures influence the behaviour of researchers, for good or for bad?
\end{enumerate}

Crucial for our current discussion is the fact that having a culture that produces a small absolute number of citations renders the discipline particularly vulnerable to manipulation of bibliometrics. With that in mind, we argue below that the answer to Q 6 is yes, and that the influence is generally not positive. Specifically, there is very good evidence that some individuals, groups, institutions and editorial boards are conspiring to tailor their publication behavior to manipulate rankings made on the basis of bibliometric analysis. 

\subsection{Bibliometrics}
Citation databases were first created as a research tool to help scholars search the academic literature by listing papers that have cited a particular, often seminal, paper. This activity has a long history dating back to the beginning of the twentieth century \cite{Godin}. However, it is generally acknowledged that the modern form of bibliometric analysis dates from the work of Garfield in the 1950s and 60s \cite{SCI-history}. The growth of this worthy endeavour into a structure whose primary purpose is the quantitative assessment of research is described in many places, among which is \cite{Pendlebury}, written by David Pendlebury, a long term employee of one of the main providers of citation data, Clarivate. Other descriptions appear in \cite{Godin} and \cite{DaMa2021}.

As Goodhart hints at in the quote at the beginning of this section, when the activity being measured involves some aspect of human behaviour, there is a large incentive to game systems. We are seeing this playing out in the field of bibliometrics across all disciplines. For example, the vulnerability of a system of research assessment based upon bibliometrics is the topic of Macdonald and Kam's paper \cite{MacdKam}. For the reasons mentioned above, the effects are particularly pernicious in the mathematical sciences.  There are ongoing efforts to develop indices that are not prone to gaming, but so far none of them have had any impact on the situation.


\subsection{The academic publishing ecosystem}
\begin{figure}[]
    \centering

\tikzset{every picture/.style={line width=0.75pt}} 

\begin{tikzpicture}[x=0.75pt,y=0.75pt,yscale=-1,xscale=1]

\draw  [fill={rgb, 255:red, 65; green, 117; blue, 5 }  ,fill opacity=0.35 ] (83,39) .. controls (83,34.58) and (86.58,31) .. (91,31) -- (182,31) .. controls (186.42,31) and (190,34.58) .. (190,39) -- (190,63) .. controls (190,67.42) and (186.42,71) .. (182,71) -- (91,71) .. controls (86.58,71) and (83,67.42) .. (83,63) -- cycle ;
\draw  [fill={rgb, 255:red, 144; green, 19; blue, 254 }  ,fill opacity=0.23 ] (19,240.8) .. controls (19,237.6) and (21.6,235) .. (24.8,235) -- (101.2,235) .. controls (104.4,235) and (107,237.6) .. (107,240.8) -- (107,258.2) .. controls (107,261.4) and (104.4,264) .. (101.2,264) -- (24.8,264) .. controls (21.6,264) and (19,261.4) .. (19,258.2) -- cycle ;
\draw  [fill={rgb, 255:red, 144; green, 19; blue, 254 }  ,fill opacity=0.23 ] (205,232.6) .. controls (205,226.75) and (209.75,222) .. (215.6,222) -- (306.4,222) .. controls (312.25,222) and (317,226.75) .. (317,232.6) -- (317,264.4) .. controls (317,270.25) and (312.25,275) .. (306.4,275) -- (215.6,275) .. controls (209.75,275) and (205,270.25) .. (205,264.4) -- cycle ;
\draw [color={rgb, 255:red, 65; green, 117; blue, 5 }  ,draw opacity=1 ][line width=1.5]    (112,85) .. controls (84.7,119.13) and (63.1,176.07) .. (51.85,212.26) ;
\draw [shift={(51,215)}, rotate = 286.99] [color={rgb, 255:red, 65; green, 117; blue, 5 }  ,draw opacity=1 ][line width=1.5]    (14.21,-4.28) .. controls (9.04,-1.82) and (4.3,-0.39) .. (0,0) .. controls (4.3,0.39) and (9.04,1.82) .. (14.21,4.28)   ;
\draw [color={rgb, 255:red, 65; green, 117; blue, 5 }  ,draw opacity=1 ][line width=1.5]    (186,82) .. controls (239.19,108.6) and (267.15,153.62) .. (277.54,205.62) ;
\draw [shift={(278,208)}, rotate = 259.32] [color={rgb, 255:red, 65; green, 117; blue, 5 }  ,draw opacity=1 ][line width=1.5]    (14.21,-4.28) .. controls (9.04,-1.82) and (4.3,-0.39) .. (0,0) .. controls (4.3,0.39) and (9.04,1.82) .. (14.21,4.28)   ;
\draw  [fill={rgb, 255:red, 74; green, 144; blue, 226 }  ,fill opacity=0.32 ] (340,37) .. controls (340,32.58) and (343.58,29) .. (348,29) -- (402,29) .. controls (406.42,29) and (410,32.58) .. (410,37) -- (410,61) .. controls (410,65.42) and (406.42,69) .. (402,69) -- (348,69) .. controls (343.58,69) and (340,65.42) .. (340,61) -- cycle ;
\draw  [fill={rgb, 255:red, 245; green, 166; blue, 35 }  ,fill opacity=0.34 ] (541,33.2) .. controls (541,28.12) and (545.12,24) .. (550.2,24) -- (633.8,24) .. controls (638.88,24) and (643,28.12) .. (643,33.2) -- (643,60.8) .. controls (643,65.88) and (638.88,70) .. (633.8,70) -- (550.2,70) .. controls (545.12,70) and (541,65.88) .. (541,60.8) -- cycle ;
\draw [color={rgb, 255:red, 65; green, 117; blue, 5 }  ,draw opacity=1 ][line width=1.5]    (201,51) .. controls (253.92,51) and (281.87,51) .. (330.03,50.06) ;
\draw [shift={(333,50)}, rotate = 178.85] [color={rgb, 255:red, 65; green, 117; blue, 5 }  ,draw opacity=1 ][line width=1.5]    (14.21,-4.28) .. controls (9.04,-1.82) and (4.3,-0.39) .. (0,0) .. controls (4.3,0.39) and (9.04,1.82) .. (14.21,4.28)   ;
\draw [color={rgb, 255:red, 163; green, 20; blue, 192 }  ,draw opacity=1 ][line width=1.5]    (198,249) .. controls (172.78,249) and (148.5,249) .. (119.69,249) ;
\draw [shift={(117,249)}, rotate = 360] [color={rgb, 255:red, 163; green, 20; blue, 192 }  ,draw opacity=1 ][line width=1.5]    (14.21,-4.28) .. controls (9.04,-1.82) and (4.3,-0.39) .. (0,0) .. controls (4.3,0.39) and (9.04,1.82) .. (14.21,4.28)   ;
\draw [line width=1.5]    (533,54) .. controls (473.9,53.02) and (472.04,53.97) .. (422.31,53.04) ;
\draw [shift={(420,53)}, rotate = 1.1] [color={rgb, 255:red, 0; green, 0; blue, 0 }  ][line width=1.5]    (14.21,-4.28) .. controls (9.04,-1.82) and (4.3,-0.39) .. (0,0) .. controls (4.3,0.39) and (9.04,1.82) .. (14.21,4.28)   ;
\draw [color={rgb, 255:red, 163; green, 20; blue, 192 }  ,draw opacity=1 ][line width=1.5]    (73,214) .. controls (81.82,159.12) and (103.13,131.13) .. (124.68,89.56) ;
\draw [shift={(126,87)}, rotate = 117.1] [color={rgb, 255:red, 163; green, 20; blue, 192 }  ,draw opacity=1 ][line width=1.5]    (14.21,-4.28) .. controls (9.04,-1.82) and (4.3,-0.39) .. (0,0) .. controls (4.3,0.39) and (9.04,1.82) .. (14.21,4.28)   ;
\draw [color={rgb, 255:red, 163; green, 20; blue, 192 }  ,draw opacity=1 ][line width=1.5]    (255,205) .. controls (232.34,135.06) and (223.27,125.28) .. (175.22,89.65) ;
\draw [shift={(173,88)}, rotate = 36.5] [color={rgb, 255:red, 163; green, 20; blue, 192 }  ,draw opacity=1 ][line width=1.5]    (14.21,-4.28) .. controls (9.04,-1.82) and (4.3,-0.39) .. (0,0) .. controls (4.3,0.39) and (9.04,1.82) .. (14.21,4.28)   ;
\draw [line width=1.5]    (402,118) .. controls (466.03,94.36) and (449.52,99.83) .. (500.62,77.06) ;
\draw [shift={(503,76)}, rotate = 156.04] [color={rgb, 255:red, 0; green, 0; blue, 0 }  ][line width=1.5]    (14.21,-4.28) .. controls (9.04,-1.82) and (4.3,-0.39) .. (0,0) .. controls (4.3,0.39) and (9.04,1.82) .. (14.21,4.28)   ;
\draw [line width=1.5]    (412,149) .. controls (456.33,119.45) and (488.04,100.57) .. (532.94,71.34) ;
\draw [shift={(535,70)}, rotate = 146.89] [color={rgb, 255:red, 0; green, 0; blue, 0 }  ][line width=1.5]    (14.21,-4.28) .. controls (9.04,-1.82) and (4.3,-0.39) .. (0,0) .. controls (4.3,0.39) and (9.04,1.82) .. (14.21,4.28)   ;
\draw [line width=1.5]    (112,230) -- (187,202) ;
\draw [line width=1.5]    (252,177) -- (262,173) ;
\draw  [fill={rgb, 255:red, 245; green, 166; blue, 35 }  ,fill opacity=0.31 ] (543,233.4) .. controls (543,227.66) and (547.66,223) .. (553.4,223) -- (619.6,223) .. controls (625.34,223) and (630,227.66) .. (630,233.4) -- (630,264.6) .. controls (630,270.34) and (625.34,275) .. (619.6,275) -- (553.4,275) .. controls (547.66,275) and (543,270.34) .. (543,264.6) -- cycle ;
\draw [color={rgb, 255:red, 206; green, 137; blue, 29 }  ,draw opacity=1 ][line width=1.5]    (429,215) .. controls (411.27,173.63) and (395.48,127.41) .. (385.45,77.29) ;
\draw [shift={(385,75)}, rotate = 78.91] [color={rgb, 255:red, 206; green, 137; blue, 29 }  ,draw opacity=1 ][line width=1.5]    (14.21,-4.28) .. controls (9.04,-1.82) and (4.3,-0.39) .. (0,0) .. controls (4.3,0.39) and (9.04,1.82) .. (14.21,4.28)   ;
\draw [color={rgb, 255:red, 197; green, 132; blue, 29 }  ,draw opacity=1 ][line width=1.5]    (62,278) .. controls (149.12,374.03) and (360.71,361.27) .. (427.03,288.23) ;
\draw [shift={(429,286)}, rotate = 130.48] [color={rgb, 255:red, 197; green, 132; blue, 29 }  ,draw opacity=1 ][line width=1.5]    (14.21,-4.28) .. controls (9.04,-1.82) and (4.3,-0.39) .. (0,0) .. controls (4.3,0.39) and (9.04,1.82) .. (14.21,4.28)   ;
\draw [color={rgb, 255:red, 197; green, 131; blue, 23 }  ,draw opacity=1 ][line width=1.5]    (252,285) .. controls (301.5,320.64) and (352.96,320.02) .. (411.23,286.04) ;
\draw [shift={(413,285)}, rotate = 149.32] [color={rgb, 255:red, 197; green, 131; blue, 23 }  ,draw opacity=1 ][line width=1.5]    (14.21,-4.28) .. controls (9.04,-1.82) and (4.3,-0.39) .. (0,0) .. controls (4.3,0.39) and (9.04,1.82) .. (14.21,4.28)   ;
\draw [color={rgb, 255:red, 0; green, 0; blue, 0 }  ,draw opacity=1 ][line width=1.5]    (579,84.57) .. controls (578.95,129.16) and (578.04,123.2) .. (578.97,201.6) ;
\draw [shift={(579,204)}, rotate = 269.31] [color={rgb, 255:red, 0; green, 0; blue, 0 }  ,draw opacity=1 ][line width=1.5]    (14.21,-4.28) .. controls (9.04,-1.82) and (4.3,-0.39) .. (0,0) .. controls (4.3,0.39) and (9.04,1.82) .. (14.21,4.28)   ;
\draw [shift={(579,81)}, rotate = 90] [color={rgb, 255:red, 0; green, 0; blue, 0 }  ,draw opacity=1 ][line width=1.5]    (14.21,-4.28) .. controls (9.04,-1.82) and (4.3,-0.39) .. (0,0) .. controls (4.3,0.39) and (9.04,1.82) .. (14.21,4.28)   ;
\draw  [fill={rgb, 255:red, 74; green, 144; blue, 226 }  ,fill opacity=0.31 ] (373,233.4) .. controls (373,227.66) and (377.66,223) .. (383.4,223) -- (460.6,223) .. controls (466.34,223) and (471,227.66) .. (471,233.4) -- (471,264.6) .. controls (471,270.34) and (466.34,275) .. (460.6,275) -- (383.4,275) .. controls (377.66,275) and (373,270.34) .. (373,264.6) -- cycle ;
\draw [line width=1.5]    (536,251) .. controls (506.6,251) and (529.06,251.96) .. (483.85,251.06) ;
\draw [shift={(481,251)}, rotate = 1.17] [color={rgb, 255:red, 0; green, 0; blue, 0 }  ][line width=1.5]    (14.21,-4.28) .. controls (9.04,-1.82) and (4.3,-0.39) .. (0,0) .. controls (4.3,0.39) and (9.04,1.82) .. (14.21,4.28)   ;
\draw [line width=1.5]    (272,168) -- (389,123) ;
\draw [line width=1.5]    (303,217) -- (400,159) ;

\draw (105,45) node [anchor=north west][inner sep=0.75pt]   [align=left] {researchers};
\draw (25,243) node [anchor=north west][inner sep=0.75pt]   [align=left] {universities};
\draw (214,234) node [anchor=north west][inner sep=0.75pt]   [align=left] {grant agencies,\\governments};
\draw (22,122) node [anchor=north west][inner sep=0.75pt]  [color={rgb, 255:red, 65; green, 117; blue, 5 }  ,opacity=1 ] [align=left] {work at};
\draw (231,92) node [anchor=north west][inner sep=0.75pt]  [color={rgb, 255:red, 65; green, 117; blue, 5 }  ,opacity=1 ] [align=left] {write \ \ \\ \ \ \ applications};
\draw (349,43) node [anchor=north west][inner sep=0.75pt]   [align=left] {journals};
\draw (549,33) node [anchor=north west][inner sep=0.75pt]   [align=left] {(commercial)\\publishers};
\draw (210,25) node [anchor=north west][inner sep=0.75pt]  [color={rgb, 255:red, 65; green, 117; blue, 5 }  ,opacity=1 ] [align=left] {publish, review,};
\draw (208,55) node [anchor=north west][inner sep=0.75pt]  [color={rgb, 255:red, 65; green, 117; blue, 5 }  ,opacity=1 ] [align=left] {serve as editors};
\draw (128,258) node [anchor=north west][inner sep=0.75pt]  [color={rgb, 255:red, 163; green, 20; blue, 192 }  ,opacity=1 ] [align=left] {evaluate};
\draw (453,30) node [anchor=north west][inner sep=0.75pt]   [align=left] {publish};
\draw (93,147) node [anchor=north west][inner sep=0.75pt]  [color={rgb, 255:red, 163; green, 20; blue, 192 }  ,opacity=1 ] [align=left] {evaluate};
\draw (187,183) node [anchor=north west][inner sep=0.75pt]  [color={rgb, 255:red, 163; green, 20; blue, 192 }  ,opacity=1 ] [align=left] {evaluate};
\draw (409.84,159.07) node [anchor=north west][inner sep=0.75pt]  [rotate=-327.28] [align=left] {pay APCs / buy bundles};
\draw (280.28,170.19) node [anchor=north west][inner sep=0.75pt]  [rotate=-338.58] [align=left] {pay APCs /bundles};
\draw (551,234) node [anchor=north west][inner sep=0.75pt]   [align=left] {indexing\\companies};
\draw (403.3,174.69) node [anchor=north west][inner sep=0.75pt]  [color={rgb, 255:red, 172; green, 115; blue, 21 }  ,opacity=1 ,rotate=-66.73] [align=left] {index};
\draw (228.13,319.8) node [anchor=north west][inner sep=0.75pt]  [color={rgb, 255:red, 172; green, 115; blue, 21 }  ,opacity=1 ,rotate=-0.87] [align=left] { buy};
\draw (585.97,108.6) node [anchor=north west][inner sep=0.75pt]  [rotate=-359.94] [align=left] {similar\\business\\goals};
\draw (381,234) node [anchor=north west][inner sep=0.75pt]   [align=left] {databases \&\\bibliometrics};
\draw (480,224) node [anchor=north west][inner sep=0.75pt]   [align=left] {produce};
\draw (316.13,287.8) node [anchor=north west][inner sep=0.75pt]  [color={rgb, 255:red, 172; green, 115; blue, 21 }  ,opacity=1 ,rotate=-0.87] [align=left] {buy};

\end{tikzpicture}    
    \caption{A simplified depiction of the academic publishing ecosystem}
    \label{fig:ecosystem}
\end{figure}
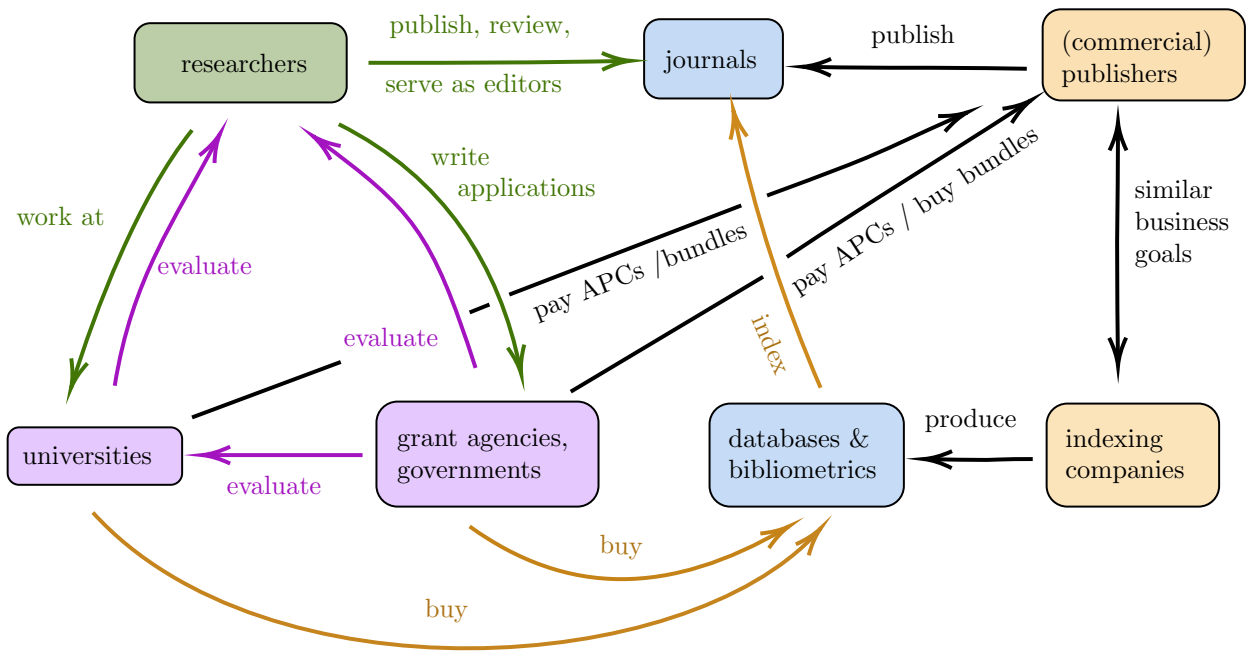

We start with a quick look at a simplified version of the academic publishing ecosystem, illustrated in Figure \ref{fig:ecosystem}. At the center of all research are, of course, the individual researchers. They work at universities or other research institutions, be they public or private. These institutions are subject to some kind of evaluation or assessment by authorities such as governments, ministries and grant agencies. Researchers publish their findings in journals, and many of them at the same time serve as editors and reviewers for these journals. Journals are published by publishers, as the name indicates. Some of these are learned societies or other not-for-profit entities but, increasingly, many of them are exceedingly profitable commercial enterprises. For simplicity, we have omitted professional societies and conferences from Figure \ref{fig:ecosystem} because they usually don't contribute to the problems that we are addressing.

Before the advent of open access (OA), these were the key players. 
University libraries would subscribe to  single journals based on the recommendations of their researchers, which limited the number of low-quality journals included in their collections. 
Digitization, OA, and the growth of the scientific community have led via a long process to different questions, depending on the perspective:
\begin{itemize}
    \item For researchers: How to find the relevant research outlet, and how to increase visibility of their work?
    \item For policymakers in governments and grant agencies: How to quantify research evaluation?
    \item For commercial publishers: How to create (more) revenue?
\end{itemize}

A first decision by large publishing houses was to promote the purchase of subscriptions in large \emph{bundles of journals} instead of individual journals---the advantage being clear, 
as long as the bundle contained the key journals of a field, the quality of a single journal mattered less, and a constant revenue was secured, as canceling single subscriptions was no longer possible. This is the policy of Springer Nature and Elsevier, who currently make up 40\% of mathematics articles indexed for example by zbMath Open \cite{COP}. OA saw the appearance of (sometimes rather substantial) \emph{article processing charges} (APCs) that authors 
pay to journals in order to publish their articles, thus creating an incentive for journals to publish as many articles as possible.

\medskip
Indexing companies, with their databases, previously existed to help researchers discover relevant papers. But the advent of quantitative assessment of research has made their databases and bibliometric products vastly more valuable. They sell these products to all larger institutions with the claim that they enable objective research assessment. 

Bibliometric analysis is possible only when a database is available that, for a sufficiently extensive set of papers ${\cal S}$, records the set $C(x)$ of papers cited by each $x \in {\cal S}$. The derivation of measures such as a journal's impact factor, an author's citation count or h-index, a field's list of HCRs or, indeed, more sophisticated statistics such as a journal ranking or an author's empirical citation distribution, is then a computational task. Since the nature of the database can have a significant effect on the measures that are computed, we shall make a few brief comments about the most important citation databases.

Clarivate \cite{clar} (known as Thomson Scientific before 2008) is a British-American publicly traded analytics company with over 10,000 employees and a revenue of US\$2.66 billion in 2022. It is the company that computes the journal impact factor and creates the list of HCRs discussed below. It offers these via the paid platform Web of Science (WoS). Every year it publishes the (paywalled) journal ranking JCR (Journal Citation Reports) based on its Science Citation Index Expanded (SCIE) database of journals and publications. The list of journals indexed is not freely available (which in itself is already problematic), but the small number of delisted journals per year underlines the fact that many predatory journals and low-quality mega-journals remain in the database \cite{Brainard}. As a rule of thumb, any journal having an impact factor is contained in WoS. 

Scopus \cite{scopus} is a product of the Dutch publishing company Elsevier. It provides journal, article and author metrics, researcher profiles, and analytical services. As the bibliometrics arm of a major publishing company, the question of whether there is a conflict of interest obviously arises. It is often claimed that Scopus' database contains questionable journals \cite{Chawla, MS}. 
It is the database used to compute the 
commercial journal ranking SJR (SCImago Journal Rank).
A quick analysis of the 2023 SJR listing in the general category ``mathematics" \cite{SJR} shows that about 20 \% of the Q1 and Q2 journals are not listed by zbMath Open---because they have little relevance for mathematics, because they are known to be of poor quality, or because they are falsely classified as being in mathematics.

Google Scholar \cite{google} is the part of Google's search engine that concentrates on the scholarly literature. It also, however, computes and publishes metrics \cite{googlemetrics}. It is often criticized for a lack of control over what is included in its database, with duplications of the same article being noticeable. Indeed at \cite{googleinclusion}, referring to a set of procedures an author can use to get their paper noticed, the phrase {\it Our search robots should normally find your paper and include it in Google Scholar within several weeks} appears. The recent study \cite{google-study} shows that the citation data of Google Scholar are easily manipulated.

zbMATH Open
is a free abstracting and reviewing service for mathematical literature and software, edited by the European Mathematical Society and FIZ Karlsruhe.
The indexing policy is to cover all available published and peer-reviewed articles, books, conference proceedings as well as other publication formats (like ArXiv preprints) pertaining to the scope defined by the MSC that present a genuinely new mathematical point of view.
Mathematical Reviews provides a subscription based service similar to zbMATH (except software and preprints) and is held in similar high regard within the mathematical science community.
In its editorial statement \cite{mathrev}, Mathematical Reviews states that it
`` is a database [the MRDB] for the mathematical sciences, produced by the American Mathematical Society, curated by mathematicians, and published on the web as MathSciNet''. 
MRDB and zbMATH Open currently cover about 1600 journals (as well as many old journals and book series) and index a total of about 5 million items.

The ArXiv \cite{arxiv} is a service provided by Cornell University with some support from other organisations. It is a repository that contains nearly 2.4 million scholarly articles in the fields of physics, mathematics, and computer science (as well as a few others). 
`Putting a paper on the ArXiv' is a first step towards publication for many mathematical scientists and a large amount of significant research can be found there.

ORCID \cite{orcid} is a global, not-for-profit organization that provides a unique, persistent identifier for individuals to use as they engage in research, scholarship, and innovation activities. Such an identifier is a useful tool for author disambiguation. ORCID has a strong set of founding principles, which we support. The relationship between ORCID and the commercial databases mentioned above is something that should be monitored as time goes by. We have already noticed instances where an ORCID identifier in, say, a login procedure has been replaced by an identifier controlled by one of the the abovementioned commercial companies.

\section{The 2023 Clarivate Exclusion}
%
\subsection{The Clarivate Exclusion}
In November 2023, Clarivate  announced that it had 
excluded the entire field of mathematics from the most recent edition of its influential list of highly-cited authors. 
Justifying this decision, it said \cite{Clarivate}
\begin{quote}
    ``(Mathematics) is a highly fractionated research domain, with few individuals working on a number of specialty topics. The average rate of publication and citation in Mathematics is relatively low, so small increases in publication and citation tend to distort the representation and analysis of the overall field. Because of this, the field of Mathematics is more vulnerable to strategies to optimize status and rewards through publication and citation manipulation.''
\end{quote}

We agree with this statement. It is interesting to see a company that supplies citation data specifically concede that such data is vulnerable to strategies to optimize status and rewards.

Some authors \cite{Catanzaro} have argued that the exclusion occurred after it became obvious that many of the people listed as HCRs in the immediately preceding years were not, in the view of the scientific community, top-level mathematicians  (see also \cite{SDC25} for further details on the exclusion). Rather they seem to be people who had manipulated their bibliographic parameters in a way that could not be ignored any more. Table \ref{table:HCR-affiliations} contains a list of the most frequent primary affiliations of the 2019 HCRs in the mathematical sciences. Note that the top listed institution, Medical University Taiwan, does \emph{not} have a program in mathematics. In the list published in November 2024, mathematics was still excluded, but without further comment. In a blog entry \cite{Clarivate2}, Head of Research Analysis David Pendlebury concedes that they excluded more than 2000 individuals for integrity concerns (the formulation is rather vague). Some remaining 6,636 individuals were listed as HCRs, but  $48\%$ of them were not associated with any definite research area (they were designated ``cross-field"). 
At least one of the people listed as cross-field in 2024, J.\,A. Tenreiro Machado, died in 2021, raising (even more) questions about the quality of the underlying data.
A recent study \cite{Chaignon} of  the set of HCRs argues that this list (as a whole, not limited to mathematics) has entered a new `phase' of development in 2019--2023: ``This phase is characterized by the constant increase in profiles
suspected of scientific misconduct, which challenges the ability of this list to identify truly influential
researchers." 

\medskip
It has been argued \cite{AET} for a long time that the commonly-used bibliometric metrics are inappropriate for the mathematical sciences, and probably for all science. Whatever your view of the reasons behind the Clarivate exclusion, its occurrence has forced the wider mathematical community to pay attention to this effect.

\raa{1.2}
\begin{table}
\begin{center}
\begin{tabular}{lc}
\toprule
University   &  \# of HCRs \\
\midrule  
 China Medical University Taiwan  & 11 \\
    King Abdulaziz 
    University, Saudi Arabia & 5 \\
    Queensland University of Technology (QUT), Australia & 3 \\ 
    Stanford University, US & 3 \\
    University of California Los Angeles, US & 3 \\ 
     Beijing Normal University, China & 2 \\ 
    Shandong Univ. of Science and Technology, China & 2 \\ 
    University of Electronic Science and Technology of China & 2 \\
    Amirkabir University of Technology, Iran & 2 \\ 
      University of Milano-Bicocca, Italy & 2 \\
    University of Urbino, Italy & 2 \\
      University of Michigan, US & 2 \\
      University of Minnesota -- Twin Cities, US & 2 \\    
       \bottomrule
       \end{tabular}
    \caption{Institutions having more than one primarily affiliated HCR in
    mathematics (from a total of 89 HCRs in math in 2019) \cite{clarHCR}}.
    \label{table:HCR-affiliations}
    \end{center}
\end{table}

\subsection{A quick look at other sciences}
%
No-one should draw the false conclusion from Clarivate's decision that mathematics has become a dubious science infiltrated by fraudsters. 
Fraudulent science occurs in other fields, not only with the aim of manipulating bibliometrics, but also in the form of claiming false results. In fact, the latter is more common in the broader sciences than in the mathematical sciences, because research papers in these disciplines often report the results of experiments, fieldwork and surveys that have to be taken on trust.
Cases of massive data falsification in the experimental sciences tend to be widely reported in the media and thus are the first that come to the public's mind.  As an example, we mention the controversy resulting in the resignation of the former president of Stanford University, neuroscientist Marc Tessier-Lavigne, in August 2023
\cite{Oransky}. 

Overall, the fact that research communities in, say,  the life sciences are much larger than in mathematics, combined with the high level of competition for substantial funding, suggests that it would not be surprising if there were a higher number of problematic publications. Several studies support this claim \cite{Byrne, Noorden2}. 
A 2020 study of retracted papers in the field of neuroscience \cite{Candal22} concluded that  ``Papers retracted originating from paper mills are increasing in frequency, posing a problem for the research community."

Systematic studies of fraudulent publishing are rare, but if one considers many events together, typical patterns and a general picture start to emerge. For example,
several scandals have shaken Spanish science in the last two years, as documented  in the international edition of the newspaper “El Pais”. The president of the old University of Salamanca, computer scientist Juan Manuel Corchado, was at the center of a very fruitful citation cartel. He now has an unprecedented 75 papers retracted \cite{Ansede3}. There were alleged paper mill purchases by
Spain’s most prolific scientific academics \cite{Ansede2}, and some universities from
Saudi Arabia paid Spanish scientists to lie about their primary affiliation
\cite{Ansede1}. One of these universities, King Abdulaziz University, is ranked second in  Table \ref{table:HCR-affiliations}. Similar scandals have been reported from other countries, including in alphabetical order:
Canada \cite{Canada1},
China \cite{China1, China2},
Hungary \cite{Gabor},
India \cite{India1},
Poland \cite{Poland1},
and Vietnam \cite{Vietnam1}.
 One of the rare exceptions of a truly systematic study is a 
detailed investigation by the German newspaper \emph{S\"uddeutsche Zeitung} and the public broadcasters NDR and WDR in 2018 \cite{SZ-Magazin}. 

\medskip
To summarize: \emph{Fraudulent publishing is a worldwide problem of substantial size in all scientific fields; it concerns researchers at all seniority levels in a variety of contexts.}


\subsection{The extremes: HCRs and retractions---and how they correlate}
\subsubsection*{A closer look at HCRs in mathematics}

In 2019, 89 researchers were listed as Highly Cited Researchers in mathematics. In 2021, Edward Dunne, then executive editor of Mathematical Reviews, conducted a thorough investigation of these researchers \cite{Dunne}, which is why our analysis will also be based on the same year. First, let us start with  a summary of Dunne's results:

\begin{enumerate}
    \item The fields in which the group of HCRs work are not representative of the mathematical sciences: for example, none is working in geometry or algebra, and only sixteen out of 60 (26\%) possible mathematical science subjects\footnote{This refers to the Mathematics Subject Classification (MSC) jointly published by Mathematical Reviews and zbMATH Open; the classes 00=General, 01=History and biography, 97=Mathematics education were left out.} actually appear as primary subject of the majority of their publications. 
    \item The most common primary subjects of HCRs are:
    partial or ordinary differential equations (37\%), statistics
    (16\%), numerical analysis (15\%), and operator theory (10\%).
    \item Instead of merely counting citations, one could use major math prizes as indicators of mathematical influence (Dunne took these to be the prizes awarded by the AMS, SIAM, EMS, and the IMU including the Fields Medal): only five of the 2019 HCRs have won any of these. \emph{The lists of HCRs and prizewinner are thus almost disjoint, despite the fact that there are 636 such winners of a total of 29 prizes, over all time}.
    \item The topics of the prizewinners are broader, they make up about 61\% of the mathematical science subjects. They are not concentrated in analysis and applications, and their most common topics (algebraic geometry and number theory) do not appear at all among the group of HCRs.
    \item The mode of `quality assessment' is totally different: HCRs are established solely on an absolute citation count of a very large, uncurated list of journals, while prizes are awarded by
    selection committees appointed by professional societies.
    \end{enumerate}

Furthermore, Dunne detected a clear difference in citation patterns. To explain it, let us introduce two similar, but crucially different numbers (see also \cite{SPA} and the critique of their work in \cite{Dunne}):

\medskip
\begin{compactdesc}
    \item[Self-citing score (SCS):] Percentage indicating the proportion of citations of an author’s works that come from the author’s own publications. The SCS is about citations {\it to} an author's papers. Specifically how many come from other papers by that author?
    
    \item[Self-referencing score (SRS):] Percentage indicating the proportion of the references in an author's paper that are to their own previous papers, as compared to referencing other authors. The SRS is about citations {\it from} an author's papers. Specifically how many are to other papers by that author?
    \end{compactdesc}

\medskip 

\raa{1.2}
\begin{table}
\begin{center}
\begin{tabular}{ccc}
\toprule
Cohort   &  median SRS & median SCS \\
\midrule  
HCRs & 15.63 & 12.91 \\
Top 1000 cited & 6.96 &6.87\\
Prizewinners & 7.91 &  0.22\\

\bottomrule
       \end{tabular}
    \caption{ Self-citing score (SCS) and self-referencing score (SRS)
    for different cohorts of mathematicians (cited from \cite{Dunne})}
    \label{table:scores}
    \end{center}
\end{table}

Care is needed for computing these numbers, and it is methodologically
reasonable and in fact necessary to restrict attention to people who won their prizes after the year 2000 (some older winners are long retired and hence do not publish any more, old citation data is often incomplete, and citation patterns could be time dependent). There is a total of 365 such mathematicians. As a third cohort for comparison, Edward Dunne and the IT team at Mathematical Reviews determined the 
1000 top-cited authors based on journal articles published between 2010 and 2020 and contained in the MR database \cite{Dunne}.

The differences are dramatic; we summarize them in Table \ref{table:scores}. While the median SRS of the top 1000 cited mathematicians is similar to that of the prizewinners, the SRS of HCRs is more than twice as large. For the median SCS, the same effect is visible and it is noticeable that prizewinners cite themselves in their subsequent  work very modestly (SCS of only 0.22).

\medskip
To summarize: \emph{HCRs cite and reference themselves about twice as often as the other two cohorts of mathematicians. Their publication data differs drastically from that of other established mathematicians.
With few exceptions, the indicator HCR is useless for detecting mathematics of good quality.}

\subsection{Retractions}
In the past, retractions were an unknown phenomenon; articles that were published were so for eternity, and mistakes were corrected in \emph{errata} published in the same journal. Then articles with falsified data began to be retracted by publishing houses to prevent the spread of incorrect data until, now, retractions are viewed as a necessary ultimate scientific correction mechanism. 

As with all fraud, only the tip of the iceberg is visible. 
As of December 2024, the database of Retraction Watch \cite{RW-data}
listed a total of 1009 articles retracted in the primary field of mathematics, and more than 3,000 who listed mathematics as a secondary topic. With a positive lens, a retraction can be a useful correction mechanism of a publishing house; as we all know, mistakes can happen. But it is fair to be suspicious when a journal or author has multiple retractions, as it could be an indicator of poor peer review  or systematic scientific misconduct. 

In the database, 
mathematics journals that have multiple retractions are mostly published by Hindawi (\emph{Journal of Function Spaces\footnote{In 2024, zbMath Open decided not to index this journal any further.}, Mathematical Problems of Engineering, Complexity}), 
which is partly due to the 2022 acquisition of Hindawi by John Wiley \& Sons and the ensuing efforts to ``clean up" the journals and increase their quality. To a far lesser extent, Sciendo Publ (\emph{Applied Mathematics and Nonlinear Sciences}) and MDPI (\emph{Symmetry, Entropy}) are affected.
Given how rare they are, it is interesting to relate the database of retractions to the list of mathematics HCRs.  Seven out of the 89 HCRs appear in the database of retractions, mainly for paper mill suspicion, plagiarism, or duplication.

Again, the problem is not limited to mathematics---and it is quickly growing.
 In 2023, more than 10,000 articles were retracted over all sciences, and (in decreasing order) the countries with most retractions per 100,000 published articles were Saudi Arabia, Pakistan, Russia, and China, and the percentage of published papers being retracted has increased from about 0.02\% in 2002 to 0.2 \% in 2023 \cite{Norden}.

\subsection{HCRs and University rankings}
In addition to serving the personal vanity of some researchers and sprucing up their CV for grant applications and salary negotiations, admission to the `Hall of Fame' of HCRs has a decisive advantage for their university: it impacts their rank  in the `Academic Ranking of World Universities' (ARWU, published since 2003), also nicknamed the `Shanghai ranking'. The computation of this rank relies on a remarkably small number of 
\href{https://www.shanghairanking.com/methodology/arwu/2024}{indicators}:
\begin{enumerate}
    \item Excellent researchers: Nobel Prizes and Fields Medals \emph{and some more}.
    \item Excellent papers: papers published in \emph{Nature} and \emph{Science}.
    \item Research output: papers indexed \emph{in major citation indices}, and the `per capita academic performance of an institution'.    
    \end{enumerate}

The second item explains part of the hype about publications in 
\emph{Nature} and \emph{Science}; for mathematics, these journals play no role.
The critical points are those emphasized in the first and third items; both rely on products of Clarivate. As the total number of Nobel Prizes and Fields Medals is rather low, the company's databases (WoS, SCIE) and bibliometrics (IF) have been used to create a new class of researchers deemed excellent, the
Highly Cited Researchers; needless to say, the major citation indices are their own. Since the database SCIE includes 
virtually all existent journals and judges them solely by their impact factor,
predatory journals and paper mills are, by construction, included as well as reputable journals.
The incentive for researchers and universities alike to manipulate 
their publication data is thus complete: a good ranking increases a university's prestige; in many countries, it is very important for students' choices of a university. It attaches a label of excellence that is rarely questioned. The problem is seriously aggravated by the fact that
it is surprisingly easy to cheat with one's affiliation, as Clarivate basically
contacts HCRs and asks for them.

And in the end, the approach is successful, at least for some. 
The top listed university `China Medical University Taiwan' from Table \ref{table:HCR-affiliations} was, in 2022, listed in the top 300 universities of the world\footnote{It ``dropped" to best 500 in 2024; one can assume that this has to do with the exclusion of mathematics, as they previously had 11 HCRs in the field, thus proving again how prone to manipulation these figures are.}. The second ranked institution welcomes 
visitors on their webpage with the opening sentence:
``King Abdulaziz University maintained its Arab leadership in international rankings, as it ranked first in the Arab world and ranked 101-150 among the best universities in the world, according to the Shanghai International Classification for the year 2021-2022, and first in the Arab world in the British QS classification index for the third year in a row." Its mathematics program is deemed to be ranked 31st in the world. The temptation to influence an institution's own position in the rankings through strategic choices is also present at prestigious traditional universities: the University of Paris-Saclay (a merger of four  `grandes \'ecoles' and  several other institutions) was founded in 2019 
with the aim of becoming a top-ranking, research-focused French university.
This, of course, works solely because the absolute size of an institution matters.

\medskip
Other rankings suffer from the same problems, be they international (like the Times Higher Education (THE) World University Ranking, which weights the indicator ``citations" by 30 \% and uses Elsevier's Scopus) or national (like NIRF in India). All are produced by private companies, based on non-curated databases riddled with mistakes and of poor quality, based on arbitrary easily-gamed criteria.

\medskip
To summarize: \emph{Rankings of universities are too easily gamed to be given any real meaning. They create false incentives and do more harm than good.}

\section{The many heads of the hydra of fraudulent publishing in mathematics}

Bibliometrics are an attempt  to quantify the quality of scientific achievements.
%
Within the mathematical sciences, it is likely that bibliometric manipulation is the most significant form of scientific fraud. Below we discuss the issues in more detail.

\medskip
Occasional dishonesty is the most widespread form of scientific misconduct. It is usually recognized only by expert peers and not approved of, but also not pursued further due to the perceived “insignificance” and difficulty of doing so. It is what we could call the \bf{``zone of occasional poor practice''\rm. Included are actions like:

\begin{itemize}
\item Papers split into least publishable units instead of a larger, more substantial article (``salami-slicing")
   \item  An exaggerated, but still (more or less) tolerable number of self-citations or citations of mathematical friends
 \item Recycling of previous text blocks, such as a standard introduction to the field (most journals now use similarity checkers and will detect these anyway)
   \item Dishonest or sloppy attribution of previous results (exaggerating one's own contribution, ``forgetting" a paper with similar content, citing a textbook instead of the original article\ldots)
   \item Sloppy  checking of existing literature
   \item Reviewers asking authors to cite the reviewers' papers
   \end{itemize}

But as D. Docampo says in \cite{Docampo}: ‘’\ldots In the complex landscape of modern academe, the maxim “publish or perish” has been gradually evolving into a different mantra: “Get cited or your career gets blighted.” (\ldots) Citation has become so important that it has driven a novel form of trickery: stealth networks designed to manipulate citations\ldots''. It is what we could call the \bf ``zone of systematic bad practice'' \rm  that includes:
\begin{itemize}
\item Citation manipulation, see for example \cite{BCLM}  where instances of manipulation through injection of meaningless texts containing a fixed set of references and addition of references during the peer-review process are documented

 \item Editors asking authors to include references to other articles in their journal after acceptance as a condition for publication   
  \item Translation and ``Copy-Paste" plagiarism
  \item Academic superiors claiming coauthorship despite not having contributed to a publication or influencing the attribution of authorship in any other way 
   \item Giving incorrect affiliations of authors on purpose\footnote{We agree that mathematical research can stretch over several affiliations. Legally, the moment of submission is the one that matters with past affiliations being acknowledged in the paper. What we mean here are such things as secondary affiliations listed as primary, and past affiliations listed so that an address looks more reputable.}
   \item Including authors without their consent\footnote{The German chemist {Karsten Krohn} from Paderborn, who died in 2013, has a paper \cite{Krohn} with nine coauthors from Pakistan in 2020, one of them wrongly claiming to be in Paderborn.}
  \item Giving incorrect funding information. While conflicts of interest because of industry-financed research is rare in mathematics, incorrect information
   is sometimes used to suggest high quality research
   \item Predatory  conferences as part of a citation cartel \cite{Nature2}
   \end{itemize}


Unfortunately, a new level of fraud can be reached when money is involved, the scientific reputation or even employment of people is directly at risk, or a substantial amount of ``criminal energy" is required for the misdeeds. 
Nick Wise, a famous science sleuth, said on the Retraction Watch Blog:
“There’s this entire economy, ecosystem of Facebook groups, Whatsapp groups, Telegram channels selling authorship for papers, selling citations, selling book chapters, selling authorship of patents” \cite{wise}.
It is what could be called the \bf ``zone of fraudulent behaviour''\rm. 
Here are some examples:

\begin{itemize}
\item Citation sales: A professional \ref{citation-broker} (with an anonymous email address) presents a list of papers that need citations and offers a 
``thank you fee"\footnote{A coworker of the first author of this article received such an offer in Summer 2024.}
\item Paper mills: It should go without saying, but buying papers from professional ghostwriters to publish under one's own name is a serious case of scientific fraud \cite{Ansede1, Ansede2, Candal22, Chambers24}
\item Authorship sales: The revenue of paper mills can be increased (or the costs for the buying person reduced) if one sells authorship on a ready-to-publish article (see for example \cite{Gabor, fbs})

\item Blackmailing: Several instances of blackmail have been documented. These include: Paper mills who threaten authors that they will reveal their scientific misconduct to journals, employers, or grant agencies
(see for example \cite{fbs})

\item 
Abuse of power: Bullying of ``non-compliant" researchers is mainly carried out by scientific seniors. 
Several cases of leaders in top academic
institutions who lost their positions for bullying have been prominent in the media in recent years

\item Identity fraud: Several cases are reported where predatory journals have included researchers as editors \emph{without their consent}, and sometimes didn't even remove their names from their webpage after complaints
\item Plagiarising articles and publishing them \emph{with an incorrect date} to make them look as if they were the original publication\footnote{Example: original article
of a group of German computer scientists from 2022 vs.\,the plagiarised article by a team of three Indian authors that appeared in a predatory journal, claimed to be from 2019, but included references published later \cite{Heider2, Heider1}.}
\item Use of pseudonyms for publications (such a case is documented in \cite{China1}), for example to avoid being connected to previous scientific misconduct
\item Non-disclosure of serious conflicts of interest
\item Creating email accounts similar to those of established researchers and suggesting them as reviewers of an article, hence being able to review one's own work.
\end{itemize}

Only rarely are fraudulent activities pursued in court. One of the few exceptions is the India based company OMICS Group Inc\footnote{We recommend reading their wikipedia article at \url{https://en.wikipedia.org/wiki/OMICS_Publishing_Group}.},
which created an intricate network of predatory journals (about 700, some of them short-lived, often rebranded) and predatory conferences (about 3000 every year). In 2016,  the U.S. Federal Trade Commission (FTC) filed a suit against OMICS. The final judgment
from May 2022 \cite{omics} found that the company had used “deceptive marketing practices” and documents many of the frauds listed above\footnote{Citations from the final judgment: ``Defendants make numerous material misrepresentations to induce consumers to submit articles to their journals"; ``Defendants misrepresent that they follow standard peer-review
practices"; ``Defendants misrepresent the impact factors of their publications"; 
``Defendants misrepresent that their publications are included in NIH’s indexing databases"; ``Defendants fail to disclose adequately their publishing fees"; ``Defendants’ deceptive conference practices" \ldots}.
Despite its name suggesting roots in the life sciences, the company also publishes several journals in mathematics.

\section{Instead of a Conclusion: A few words about AI}

 While AI can be an invaluable tool for a variety of tasks, it also introduces new
 risks and challenges that the research community is only beginning to comprehend. 

\medskip

Most journal editors and publishers concur that authors may use AI-driven tools such as DeepL, ChatGPT, or DeepSeek for copy-editing and proofreading purposes. This practice can be particularly beneficial for non-native English speakers. However, it is critical that authors take responsibility for ensuring that the content of their publications (including citations, tables and figures) is not altered by automated translation or editing. Under normal circumstances, the use of these AI services for copyediting does not need to be formally acknowledged in the publication.
In contrast, when AI is used for more substantive tasks--such as generating or interpreting research content--many publishers now require an explicit statement disclosing this use, and no AI may be listed as a co-author.
This policy shift underscores a growing concern that AI-generated text could be used to fabricate or inappropriately alter scientific claims, thereby undermining scientific integrity. The appendix to \cite{Glynn} provides a valuable
compilation of excerpts from AI policy statements by publishers and scientific organizations worldwide, demonstrating a certain convergence on the key principles outlined above.

\medskip

The advent of sophisticated AI language models has made it cheaper and easier to produce fake research. Paper mills are likely to exploit AI to generate plausible-sounding texts with minimal human oversight \cite{ThePublicationPlan2024, Nature2023}. 
%
%
It also presents a major challenge for fraud detection: As AI-generated texts become more coherent and context-aware, it becomes increasingly difficult for traditional screening methods to distinguish them from human work. This has spurred investment in new “integrity software”—AI tools being developed to catch signs of automated content generation. However, this will almost certainly lead to an arms race: as detection tools improve, so do the generative models, potentially creating a never-ending cycle of innovation on both sides.
At present, \ref{tortured} (unusual or contorted language constructions) remain a common red flag which can be used to identify AI-generated text \cite{Cabanac} (see Table \ref{table:tortured}). Yet, relying on textual anomalies alone may not be sustainable as AI models become more refined. The community is thus calling for robust guidelines, shared databases of fraudulent works, and transparent protocols for both detection tools and editorial processes. 

\medskip

A significant turning point was Wiley’s experience with Hindawi, which reportedly cost the company substantial resources in handling fraudulent submissions \cite{RW-Wiley}. Consequently, Wiley announced a pilot “AI-powered Paper Mill Detection” service \cite{WileyNewsRelease2024}. Other major publishers, such as Springer, are following suit, looking for comprehensive strategies that combine algorithmic checks, expert review, and community-based reporting \cite{COPE}.
In his recent investigation of undeclared AI usage in scientific publications \cite{Glynn}, the author Alex Glynn concludes: ``This analysis [\ldots] reveals that the problem is widespread, penetrating the journals and conference proceedings of
highly respected publishers".

\medskip

Beyond publishers, funding agencies have also begun to issue guidelines on AI. The U.S. National Science Foundation (NSF), for instance, has released statements emphasizing the importance of transparency in AI usage and the need for responsible use, including the prohibition of  ``uploading any content from proposals, review information and related records to non-approved generative AI tools" \cite{NSF2024}.
This raises the even more fundamental question of what it means for the mathematical community if LLMs are trained with our papers, with or without our consent \cite{Harris}.

\medskip

\emph{To summarize: AI will exacerbate the existing problems of scientific publishing, and the options for cheating will multiply. It remains to be seen whether the scientific community has the resources necessary to counteract this.  
In the meantime, we appeal to everyone to stand up for the strengthening of ethical standards in their fields and to support honest scientists in this unequal struggle.}

\appendix\section{Glossary}

{     
\raa{1.2}
\begin{table}
\begin{center}
\begin{tabular}{cc}
\toprule
Tortured phrases & original meaning\\
\midrule
halfway differential condition & partial differential equation\\
worldwide constant & global constant\\
phantom examinations & spectral analysis\\
enormous information	& big data \\
fluffy rationale & fuzzy logic\\
hereditary calculation & genetic algorithm\\
man-made intelligence &	artificial intelligence\\
mistaken back spread & error back propagation\\
strategic relapse& logistic regression\\
warmth transmission	& heat transfer\\
beginning boundaries &	initial parameters\\
likelihood dissemination & probability distribution\\
informational collection & data set\\
distant detecting &	remote sensing\\
quick Fourier transform	& Fast Fourier Transform\\
JPEG pressure  & JPEG compression\\
pre- and post-handling & 	pre- and post-processing\\
\bottomrule
\end{tabular}
\caption{Collection of \ref{tortured}}\label{table:tortured}
\end{center}
\end{table}
}
    
\begin{compactdesc}
\item[\namedlabel{citation-cartel}{Citation cartel}.] Group of individuals that agree to cite authors from the cartel and their preferred journals, regardless of their relevance for their work. 
\item[\namedlabel{citation-broker}{Citation broker}.] Anonymous e-mailer who offers a``thank you fee" for each citation of a given list of articles. Payment is typically in cryptocurrency.
\item[\namedlabel{citejacked}{Citejacked journal}.] 
Established journal citing articles in a \ref{hijacked-v1}.
\item[\namedlabel{mega-journal}{Mega-journal}.] Journals publishing huge numbers of articles per year; they did not exist before the advent of Open Access. One of the most prolific mega-publishers is MDPI, founded in Basel, Switzerland. It currently owns 433 journals. Its top publication \emph{International Journal of Environmental Research and Public Health} publishes nearly 17 000 studies each year (it was delisted from WoS in 2023 and will hence lose its impact factor \cite{Brainard}). Its largest math journal \emph{Mathematics} has more than 1000 people in the editorial board (27 of them are HCRs) and publishes over 6 000 articles in 2023 --- compared to a total of about 120 000 papers reviewed in zbMath in the same year. Its yearly revenue is estimated at approximately 10 million SF.
Multiple scandals and retractions prove that it is just not possible to maintain high scientific quality in this huge quantity of publications (see for example \cite{RW-MDPI} or \cite{BishopBlog}; the last one illustrates again how articles of bad quality manage to make it into newspapers).
\item[\namedlabel{hijacked-v1}{Hijacked journal, v1}.] A fake website that uses the titles, ISSNs, and metadata of existing established journals. These first emerged around 2010; they are closely related to \ref{citejacked}. So far, about 250 cases are documented, none from mathematics; recently, journals from Springer Nature and Elsevier have been targeted \cite{hijacked}. In a recent study, about 66 \% of articles in hijacked journals contained evidence of
plagiarism \cite{Albakina}.

\item[\namedlabel{hijacked-v2}{Hijacked journal, v2}.] A journal whose editorial board has been overtaken or infiltrated by members of a citation cartel. Sometimes, the journal becomes permanently predatory.
\item[Paper mill.]
A commercial organisation that engages in the large scale production and sale of papers to researchers, academics, and students 
for publication in peer reviewed journals. Many paper mill papers include fabricated data
\cite{Candal22}. A detailed analysis of a Russian paper mill is given in \cite{Albakina23}.
\item[Predatory journal.] 
Journal exploiting the open-access publication model to deceive authors into paying them a fee. These publishers often lie about the journal’s impact factor, have poor editorial standards and falsely claim to provide a rigourous peer-review process \cite{Nature}.

A famous example is the journal \href{https://zbmath.org/serials/?q=Advances+in+Difference+Equations}{Advances in Difference Equations}; zbMath statistics clearly reveal the presence of a strong citation cartel. After several failed attempts to renew it, SpringerNature decided to close it. The last paper was published in 2021.
There is no list of predatory journals; however, recommendations on this issue may be found in \cite{WG-2}.

\item[Predatory conference.]
A conference with weak or no peer review for presentations, poor organization and a focus on making money for the organizers. They are typically taking place in large conference venues, with no or very few people actually present. Emphasis is often on virtual participation options and the promise of subsequent publication indexed in reputable citation databases.
\item[\namedlabel{sleuth}{Scientific Sleuth}.]
Persons looking for problems in the scientific literature, often in their spare time
\cite{RW-sleuths}.
\item[\namedlabel{tortured}{Tortured phrases}.]
Incorrect scientific vocabulary that stems from machine translation / paraphrase / generation of texts \cite{Cabanac, CLM21}. They are the indicator used by the 
Problematic Paper Screener \cite{PPS}. A collection of tortured phrases from mathematics and data science (collected from real papers flagged by the Problematic Paper Screener) may be found in Table \ref{table:tortured}.

\end{compactdesc}


\end{document}